\documentclass[11pt]{article}

\usepackage[bottom]{footmisc}

\usepackage[latin1]{inputenc}
\usepackage[all]{xypic}

\usepackage{amsmath,amssymb,amsfonts,enumerate,bbm,eucal}

\newcommand{\tmop}[1]{\operatorname{#1}}

\newcommand{\toisom}{\stackrel{\sim}{\to}}

\DeclareMathOperator*{\liminv}{\underleftarrow{\lim}}

\newcommand{\tensor}{\otimes}

\newcommand{\ideal}[1]{{\mathfrak #1}}

\newcommand{\spec}[1]{{\tmop{Spec}\,({#1})}}

\newcommand{\id}{{\tmop{id}}}

\newenvironment{enumerateroman}{\begin{enumerate}[i)]}{\end{enumerate}}

\DeclareMathOperator{\Tor}{Tor}

\DeclareMathOperator{\im}{im}

\newtheorem{theorem}{Theorem}[section]
\newtheorem{lemma}{Lemma}[section]

\author{J\"urgen B\"ohm}
\title{A new proof of the local criterion of flatness}

\begin{document}

\date{\today}

\maketitle

\begin{abstract}
Let $(A,\ideal{m}_A) \to (B,\ideal{m}_B)$ be a local morphism of local
noetherian rings and $M$ a finitely generated $B$--module. 
Then it follows from $\Tor^A_1(M,A/\ideal{m}_A) = 0$ that
$M$ is a flat $A$--module. This is usually called the 
"local criterion of flatness". We give a proof that proceeds along
different lines than the usual textbook proofs, 
using completions and only elementary properties of
flat modules and the $\Tor$--functor.
\end{abstract}

\section{A new proof of the local criterion of flatness}

\subsection{Introduction}

It is well known, that for a finite $A$--module $M$ over a noetherian local ring $A$
the truth of $\Tor^A_1(M, A/\ideal{m}_A) = 0$ implies that $M$ is $A$--flat.

Under some assumptions on $M$ this can be generalized to the case where $M$ is no longer
finitely--generated over $A$. The respective theorems are often called 
"local criterion of flatness" (see for example \cite[(20.C) Theorem 49]{matsumura:70},
\cite[Theorem 6.8]{eisenbud:95})

The proof here given proceeds along different lines than the proofs cited above, making
essential use of completions and only the simplest properties of the $\Tor$--functor
and of flat modules.

For all results in commutative algebra we refer to
\cite{eisenbud:95}, \cite{matsumura:70} and \cite{atiyah:69}.

\subsection{An introductory lemma}

\begin{lemma}
\label{lem:artinian-projectivity}
Let $(A,\ideal{m}_A)$ be an artinian local ring and $M$ an $A$--module. Then the following assertions are equivalent
\begin{enumerateroman}
\item
$M$ is free.
\item
$M$ is projective.
\item
$M$ is $A$-flat.
\item
$\Tor^A_1(M,k_A) = 0$, where $k_A = A/\ideal{m}_A$.
\end{enumerateroman}
\end{lemma}

{\bf Proof.}
We prove iv) $\Rightarrow$ i). First by considering the sequence 
$$0 \to \ideal{m}^p/\ideal{m}^{p+1} \to A/\ideal{m}^{p+1} \to A/\ideal{m}^p \to 0$$
one concludes that
$$\Tor^A_1(M, A/\ideal{m}^p ) = 0$$
for all $p \geqslant 0$.

Now assume we have found inductively an isomorphism:
$$\gamma: F \tensor_A A/\ideal{m}^p \toisom M \tensor_A A/\ideal{m}^p$$
where $F = \bigoplus_{i \in I} A$ is a free $A$--module.

First construct a small commutative diagram of $A$--modules
\begin{equation}
\xymatrix@-0.4pc{
 & & 0 \\
0 \ar[r] & F \tensor_A A/\ideal{m}^p \ar[r]^{\gamma} & M \tensor_A A/\ideal{m}^p \ar[r] \ar[u] & 0 \\
 & F \ar[u] \ar[r]^{\zeta} & M \tensor_A A/\ideal{m}^{p+1} \ar[u]
}
\end{equation}

By successive tensoring we can construct a diagram:
\begin{equation}
\xymatrix@-0.4pc{
      & & 0 & 0 & 0 \\
0 \ar[r] & 0 \ar[r] & F \tensor_A A/\ideal{m}^p \ar[r]^{\gamma} \ar[u] &
				M \tensor_A A/\ideal{m}^p \ar[r] \ar[u] & 0 \ar[r] \ar[u] & 0 \\
0 \ar[r] & L \ar[r] \ar[u] & F \tensor_A A/\ideal{m}^{p+1} \ar[r]^{\delta} \ar[u]^{\beta} &
				M \tensor_A A/\ideal{m}^{p+1} \ar[r] \ar[u]^{\alpha} & P \ar[r] \ar[u] & 0 \\
0 \ar[r] & 0 \ar[r] \ar[u] & F \tensor_A \ideal{m}^p/\ideal{m}^{p+1} \ar[r]^{\chi} \ar[u]^{\phi} &
				M \tensor_A \ideal{m}^p/\ideal{m}^{p+1} \ar[r] \ar[u]^{\psi} & 0 \ar[r] \ar[u] & 0 \\
      & 0 \ar[u] & 0 \ar[u] & 0 \ar[u] & 
}
\end{equation}

Note, especially, that $\chi = \zeta \tensor_A \id_{k_A} \tensor_{k_A} \id_{\ideal{m}^p/\ideal{m}^{p+1}}$
is an isomorphism, because $\zeta \tensor_A \id_{k_A}$ is one.

Now by the snake lemma, for example, follows $L = 0$ and $P = 0$. 
So we can conclude from the fact that $\gamma$ is an isomorphism, that $\delta$ is an isomorphism too.

So, by induction, 
$M \tensor_A A/\ideal{m}^p$ is a free $A/\ideal{m}^p$--module
for all $p \geqslant 0$. As $A$ is artinian and $\ideal{m}^{p_0} = 0$ for a certain $p_0 > 0$, we have
that $M$ is $A$--free.

\subsection{The main theorem}

\begin{theorem}
\label{thm:main-theorem}
Let $(A,\ideal{m}_A) \to (B,\ideal{m}_B)$ be a local morphism of noetherian local rings.
Further let $Y$ be a finitely generated $B$--module.
Then the following assertions are equivalent:

\begin{enumerateroman}
\item
$Y$ is a flat $A$--module.
\item
$\Tor^A_1(Y, k_A) = 0$, where $k_A = A/\ideal{m}_A$.
\end{enumerateroman}
\end{theorem}

{\bf Proof.}
We prove the nontrivial direction:
Consider an exact sequence of $A$--modules
\begin{equation}
\label{eq:base-seq}
0 \to N \to F \to Y \to 0
\end{equation}
where $F = \bigoplus_{i\in I} A$ is a free $A$--module.

Now consider in \eqref{eq:base-seq} the filtrations
\begin{align}
N_k = N \cap \ideal{m}_A^k F, & & \ideal{m}_A^k F = F_k, & & Y_k = \ideal{m}_A^k Y
\end{align}
From them result exact sequences of $\widehat{A}$--modules, where $\widehat{A}$ is the
$\ideal{m}_A$--adic completion of $A$:
\begin{equation}
\label{eq:compl-seq-1}
\xymatrix{
0 \ar[r] & \widehat{N} \ar[r] \ar[d] & \widehat{F} \ar[r] \ar[d] & \widehat{Y} \ar[r] \ar[d] &
 0 \\
0 \ar[r] & N/N_k \ar[r] & F/F_k \ar[r] & Y \tensor_A A/\ideal{m}_A^k \ar[r] & 0 
}
\end{equation}
There the lower sequences are exact by definition of $N_k$, $F_k$ and $Y_k$.

\paragraph*{An artinian interlude}
We use the following Lemma
\begin{lemma}
\label{lem:artinian-lemma}
Under the conditions of the the theorem it follows from $\Tor^A_1(Y,k_A)=0$ that 
$Y/\ideal{m}_A^k Y$ is a projective $A/\ideal{m}_A^k$--module.
\end{lemma}
Call $A' = A/\ideal{m}_A^k$. Then by lemma \ref{lem:artinian-projectivity} it is enough to prove, that 
$$\Tor^{A'}_1(Y/\ideal{m}_A^k Y, A'/\ideal{m}_{A'}) = 0.$$

Now from the sequences 
$$0 \to \ideal{m}_A^p/\ideal{m}_A^{p+1} \to A/\ideal{m}_A^{p+1} \to A/\ideal{m}_A^p \to 0$$
and the base assertion $\Tor^A_1(Y,A/\ideal{m}_A) = 0$ it follows inductively, that
$$\Tor^A_1(Y,A/\ideal{m}_A^{p}) = 0$$
for all $p \geqslant 0$.

Now consider the exact sequence 
$$0 \to L \to F \to Y \to 0$$
with free $A$--module $F$. Tensoring with $A/\ideal{m}_A^k$ gives the exact sequence
\begin{equation}
\label{eq:test-seq}
0 \to L \tensor_A A/\ideal{m}_A^k \to F \tensor_A A/\ideal{m}_A^k \to Y \tensor_A A/\ideal{m}_A^k \to 0.
\end{equation}
As $F \tensor_A A/\ideal{m}_A^k$ is a free $A/\ideal{m}_A^k$ module it follows that
$\Tor^{A'}_1(Y/\ideal{m}_A^k Y, k_{A'}) = 0$ if tensoring \eqref{eq:test-seq} by $- \tensor_{A'} k_A$ results
in an exact sequence. But \eqref{eq:test-seq} $\tensor_{A'} k_A$ is nothing else but
$$0 \to L \tensor_A k_A \to F \tensor_A k_A \to Y \tensor_A k_A \to 0$$
which is exact by $\Tor^A_1(Y,k_A) = 0$.

\paragraph*{Climbing the ladder}
We can therefore find in the lower sequences of 
\eqref{eq:compl-seq-1} a projective splitting of $F/F_k$ that climbs from 
$k$ to $k+1$:
\begin{equation}
\label{eq:climb-ladder}
\xymatrix{
0 \ar[r] & N/N_{k+1} \ar[r]^{\phi_{k+1}} \ar[d]^{\alpha} & 
	F/F_{k+1} \ar[r]^{\psi_{k+1}} \ar[d]^{\beta} &
	Y/\ideal{m}_A^{k+1}Y \ar[r] \ar[d]^{\gamma} \ar@/^/[l]^{s_{k+1}}& 0 \\
0 \ar[r] & N/N_{k} \ar[r]^{\phi_{k}} & 
	F/F_{k} \ar[r]^{\psi_{k}} &
	Y/\ideal{m}_A^{k}Y \ar[r] \ar@/^/[l]^{s_{k}} & 0 \\	
}
\end{equation}
First construct $s_{k+1}'$ from the condition $\psi_{k+1} s_{k+1}' = \id$. From this
follows $\psi_k \, (\beta s_{k+1}' - s_{k} \gamma) = 0$. So we have a map
$s_{k+1}'' = \beta s_{k+1}'-s_{k} \gamma:Y/Y_{k+1} \to N/N_k$. We lift it to a map
$s_{k+1}''':Y/Y_{k+1} \to N/N_{k+1}$. Then $s_{k+1} = s_{k+1}' - \phi_{k+1} s_{k+1}'''$.
is a lifting of $s_k$ that makes the diagram \eqref{eq:climb-ladder} commute.

\paragraph*{The splitting diagram}
So we get from above a commutative diagram with exact rows:
\begin{equation}
\label{eq:split-diagram}
\xymatrix{
0 \ar[r] & N/N_{k} \ar[r]^{\phi_{k}} & 
	F/F_{k} \ar[r]^{\psi_{k}} &
	Y/\ideal{m}_A^{k}Y \ar[r] \ar@/^/[l]^{s_{k}} & 0 \\	
0 \ar[r] & \widehat{N} \ar[r] \ar[u] &
 \widehat{F} \ar[r] \ar@/^/[l]^{q} \ar[u] &
  \widehat{Y} \ar[r] \ar@/^/[l]^{s} \ar[u] & 0
}
\end{equation}

\paragraph*{Using the splitting}

Let $M$ be a finitely generated $A$--module and $\widehat{M}$ its completion. We write
$$\widehat{M} = \liminv_k M_{(k)} = \liminv_k M/\ideal{m}_A^k M.$$

Now consider the mapping
\begin{equation}
\xymatrix@-0.4pc{
(\liminv_k F_{(k)}) \tensor_{A} M \ar[r]^{u_M} \ar[d] &
\liminv_k (F_{(k)} \tensor_{A} M) \ar[d] \\
F_{(k)} \tensor_A M \ar[r] & (F_{(k)} \tensor_A M)
}
\end{equation} 
where we made use of the abbreviation $X_{(k)} = X/\ideal{m}_A^k X$.

Note also
\begin{equation}
\liminv_k (F_{(k)} \tensor_{A_{(k)}} M_{(k)})
 = \liminv_k (F \tensor_A M_{(k)}) = \liminv_k (F_{(k)} \tensor_A M)
\end{equation}

We will prove that $u_M$ is an isomorphism for a finitely generated $A$--module $M$.

First we prove this for $M = E$ free of rank $r$:
\begin{equation}
\liminv_k F_{(k)} \tensor_{A} A^r = (\liminv_k F_{(k)})^r = 
\liminv_k F_{(k)}^r = \liminv (F_{(k)} \tensor_{A} A^r)
\end{equation}

Now consider a presentation
$$E' \to E \to M \to 0$$
with finite rank free $E$, $E'$ and the diagram
\begin{equation}
\label{eq:diagram-isom}
\xymatrix@-0.4pc{
(\liminv F_{(k)}) \tensor_{A} E'\ar[r] \ar[d]^{\sim} &
(\liminv F_{(k)}) \tensor_{A} E \ar[r] \ar[d]^{\sim} &
(\liminv F_{(k)}) \tensor_{A} M \ar[r] \ar[d]^{u_M} & 0 \\
\liminv (F \tensor_A E'_{(k)}) \ar[r]  &
\liminv (F \tensor_A E_{(k)}) \ar[r]  &
\liminv (F \tensor_A M_{(k)}) \ar[r] & 0
}
\end{equation}

It proves that $u_M$ is an isomorphism, if we can show, that the bottom row is exact.
We will show this in a moment, but first we will further the main line of argument:

Consider the line
\begin{multline}
\label{eq:tensor-F-einbett}
\widehat{F} \tensor_{A} M = 
(\liminv_k F_{(k)}) \tensor_{A} M = \liminv_k (F_{k} \tensor_{A} M_{(k)}) =
\liminv_k \bigoplus_i M_{(k)} \hookrightarrow \\
 \hookrightarrow \liminv_k \prod_i M_{(k)} = \prod_i \liminv_k M_{(k)} =
\prod_i \widehat{M}
\end{multline}

So we have a canonical injection
\begin{equation}
\label{eq:tensor-F-einbett-1}
\widehat{F} \tensor_A M \hookrightarrow \prod_i \widehat{M}.
\end{equation}

In the following we use the so called "Mittag--Leffler property" of inverse systems.
See for example \cite[Proposition II.9.1.]{hartshorne:83} for definition of 
and elementary facts about this property.

It remains to prove the exactness of the lower row in diagram \eqref{eq:diagram-isom}:

Start with the exact sequences
\begin{equation}
E'_{(k)} \to E_{(k)} \to M_{(k)} \to 0 
\end{equation}
which form an inverse system in $k$. Splice them into short exact sequences
\begin{align}
\label{eq:ML-sequences-1}
0 & \to P_k \to E'_{(k)} \to Q_k \to 0 \\
\label{eq:ML-sequences-2}
0 & \to Q_k \to E_{(k)} \to M_{(k)} \to 0
\end{align}
The above two systems of sequences each form an inverse system. We note that 
{\em $(P_k)_k$ and $(Q_k)_k$ have the Mittag--Leffler property (ML) as they consist of Artin--modules
only}.

Now tensoring the sequences in \eqref{eq:ML-sequences-1}, \eqref{eq:ML-sequences-2}
with $\tensor_A F$ retains their exactness
{\em and the (ML)--property on $(P_k \tensor_A F)_k$ and $(Q_k \tensor_A F)_k$ too}.

This is because, if $\psi_{k'k}: P_{k'} \to P_k$ and we have $\im (\psi_{k'k}) = \im (\psi_{k''k})$, then
$\im (\psi_{k'k}\tensor \id_F ) = \im (\psi_{k''k} \tensor \id_F)$ too.

Now take a sequence $0 \to M' \to M$ of two finitely generated $A$--modules.

Then consider the diagram
\begin{equation}
\xymatrix@-0.4pc{
 &  & 0 \ar[d] & 0 \ar[d] \\
0 \ar[r] & X \ar[r] & M' \tensor_A \widehat{F} \ar[d] \ar[r] & 
	M \tensor_A \widehat{F} \ar[d] \\
& 0 \ar[r]  & \prod_{i \in I} \widehat{M}' \ar[r]
		& \prod_{i \in I} \widehat{M}
}
\end{equation}
From this follows $X=0$ and therefore the conclusion that
tensoring with $\tensor_A \widehat{F}$ is exact on injections of finitely generated 
$A$--modules. 

From this follows at once that $\widehat{F}$ is
a flat $A$--module.

 As $\widehat{Y}$ and $\widehat{N}$ are split-summands of $\widehat{F}$ they
are $A$--flat too:

\begin{lemma}
The modules $\widehat{F}$, $\widehat{Y}$, $\widehat{N}$ from above are $A$--flat modules.
\end{lemma}

Additionally we have
\begin{lemma}
\label{lem:faithfully-flat}
The module $\widehat{B}$ is a faithfully--flat $B$--module. There $\widehat{B}$ is the
completion of $B$ with respect to the filtration $(\ideal{m}_A^k B)$, that is
$\liminv_k B/\ideal{m}_A^k B$.
\end{lemma}

First $\widehat{B}$ is $B$--flat as the $\ideal{m}_A B$--adic completion of $B$. Furthermore, we have
$\ideal{m}_B \supset \ideal{m}_A B$ and therefore 
$\widehat{\ideal{m}}_B = \ideal{m}_B \widehat{B}$ is a maximal ideal in $\widehat{B}$. It has the
property, that under $B \to \widehat{B}$ we have $\widehat{\ideal{m}}_B \cap B = \ideal{m}_B$.

So, together with the going--down property for flat extensions, we conclude, that
$\spec{\widehat{B}} \to \spec{B}$ is surjective and therefore $\widehat{B}$ a 
faithfully $B$--flat module.

\begin{lemma}
\label{lem:yhat-as-tensor}
There is an isomorphism $\widehat{Y} = Y \tensor_B \widehat{B}$
\end{lemma}

This is well known (\cite[Proposition 10.13]{atiyah:69}).

\paragraph*{The conclusion}

Now consider an injection of two finitely generated 
$A$--Modules $0 \to N' \to N$. Tensoring with $Y$ gives the exact sequence
\begin{equation}
0 \to P \to N' \tensor_A Y \to N \tensor_A Y
\end{equation}
Tensoring with $\tensor_B \widehat{B}$ leads to
\begin{equation}
0 \to P \tensor_B \widehat{B} \to N' \tensor_A (Y \tensor_B \widehat{B})  \to 
			N \tensor_A (Y \tensor_B \widehat{B})
\end{equation}
As $Y \tensor_B \widehat{B}$ is $\widehat{Y}$ and $\widehat{Y}$ is $A$--flat it 
follows, that $P \tensor_B \widehat{B} = 0$.

Now by lemma \ref{lem:faithfully-flat} we conclude $P = 0$. So $Y$ is a flat $A$--module.


\bibliography{mainlitbank}

\begin{thebibliography}{1}

\bibitem{atiyah:69}
Michael~F. Atiyah and I.G. Macdonald.
\newblock {\em {Introduction to commutative algebra.}}
\newblock {Reading, Mass.-Menlo Park, Calif.- London-Don Mills , Ont.: Addison-
  Wesley Publishing Company }, 1969.

\bibitem{eisenbud:95}
David Eisenbud.
\newblock {\em {Commutative algebra. With a view toward algebraic geometry.}}
\newblock {Graduate Texts in Mathematics. 150. Berlin: Springer-Verlag}, 1995.

\bibitem{hartshorne:83}
Robin Hartshorne.
\newblock {\em {Algebraic geometry. Corr. 3rd printing.}}
\newblock {Graduate Texts in Mathematics, 52. New York-Heidelberg-Berlin:
  Springer- Verlag. XVI, 496 p.}, 1983.

\bibitem{matsumura:70}
H.~Matsumura.
\newblock {\em {Commutative algebra.}}
\newblock {Mathematics Lecture Note Series. New York: W. A. Benjamin, Inc. xii,
  262 p. }, 1970.

\end{thebibliography}

\end{document}